\documentstyle[12pt,epsf]{article}

\advance \topmargin by -\headheight
\advance \topmargin by -\headsep

\evensidemargin \oddsidemargin
 \def\ex{{\hbox{\rm e}}}

  \def\tr{{\hbox{\rm Tr}}}

 \def\vev{vacuum expectation value}

\def\knot{Journal of Knot Theory and Its Ramifications}

\newcommand{\QQ}{{\mbox{{\bf Q}}}}
\newcommand{\RR}{{\mbox{{\bf R}}}}
\newcommand{\RRs}{{\mbox{{\small \bf R}}}}

\newcommand{\W}{{\cal W}}

\newfont{\namefont}{cmr10}
\newfont{\addfont}{cmti7 scaled 1440}
\newfont{\headfontb}{cmbx10 scaled 1728}
\begin{document}
\begin{titlepage}
\begin{center}
{\headfontb Primitive Vassiliev Invariants and Factorization} 
\vskip0.2cm
{\headfontb in Chern-Simons Perturbation Theory}
\footnote{This work is supported in part by funds provided by the
U.S.A. DOE under cooperative research agreement
\#DE-FC02-94ER40818 and by the DGICYT of Spain under grant PB93-0344.}
\end{center}
\vskip 0.3truein
\begin{center}
{\namefont M. Alvarez\footnote{Email: marcos@mitlns.mit.edu}}
\end{center}
\begin{center}
{\addfont{Center for Theoretical Physics,}}\\
{\addfont{Massachusetts Institute of Technology}}\\
{\addfont{Cambridge, Massachusetts 02139 U.S.A.}}
\end{center}
\vskip 0.3truein
\begin{center}
{\namefont J.M.F. Labastida\footnote{Email: labastida@gaes.usc.es}}
\end{center}
\begin{center}
{\addfont{Departamento de F\'\i sica de Part\'\i culas,}}\\
{\addfont{Universidade de Santiago}}\\
{\addfont{E-15706 Santiago de Compostela, Spain}}
\end{center}
\vskip 1truein

\begin{center}
\bf ABSTRACT
\end{center}
The general structure of the perturbative expansion of the vacuum expectation 
value of a Wilson line operator in Chern-Simons gauge field theory is 
analyzed. The expansion is organized according to the independent group 
structures that appear at each order. It is shown that the analysis is greatly 
simplified if the group factors are chosen in a certain way that we call 
canonical. This enables us to show that the logarithm of a polinomial knot 
invariant can be written in terms of primitive Vassiliev invariants only.  
\vskip3.5truecm
\leftline{MIT-CTP-2524}
\leftline{USC-FT/16-96}
\leftline{q-alg/9604010   \hfill April 1996}
\smallskip
\end{titlepage}
\setcounter{footnote}{0}


\section{Introduction} 

Vassiliev invariants, or numerical invariants of
finite type, are a set of knot invariants first proposed  in \cite{vassi}
to classify knot types.  To each knot corresponds an infinite sequence of 
rational numbers which have to satisfy some consistency conditions in order 
to be knot class invariants. This infinite sequence is divided into finite 
subsequences, each one characterized by its order, which form vector spaces. 
The number of independent elements in each finite 
subsequence is called the dimension of the space of Vassiliev invariants at 
that order. 

Apart from the original definition \cite{vassi} of these invariants, there are 
other approaches to the subject. Since their formulation in terms of inductive 
relations for singular knots \cite{birlin,birman} and of
their relation to knot invariants based on quantum groups or in
Chern-Simons gauge theory \cite{drorcon,drortesis,birlin,lin,birman}, several 
works have been
performed to analyze Vassiliev invariants in both frameworks 
\cite{konse,drortopo,piuni,numbers,haya}. In \cite{konse,drortopo} it was shown
that Vassiliev invariants can be understood in terms of representations of
chord diagrams whithout isolated chords modulo the so called 4T relations
(weight systems), and that using semi-simple Lie algebras  weight systems can
be constructed. It was also shown in \cite{drortopo}, using Kontsevitch's
representation for Vassiliev invariants \cite{konsedos} that the space of
weight systems is the same as the space of Vassiliev invariants. In 
\cite{piuni} it was argued that these representations are precisely the ones 
underlying quantum-group or Chern-Simons invariants. 

 We observed in \cite{numbers} that the generalization of 
the integral or geometrical knot invariant first proposed in \cite{gmm} and 
further analyzed in \cite{drortesis}, as well as the invariant itself, are 
Vassiliev invariants. In \cite{numbers} we proposed an organization of those
geometrical invariants and we described a procedure for their calculation
from known polynomial knot invariants. This procedure was applied to
obtaining Vassiliev knot invariants up to order six for all prime knots up
to six crossings. These geometrical invariants have also been studied by
Bott and Taubes \cite{bot} using a different approach. The relation of this
approach to the one in \cite{numbers} has been studied recently 
in \cite{altfr}.

An interesting outcome of the analysis presented in \cite{numbers} is
the well known fact that 
the Vassiliev invariants of a given knot 
form an algebra in the sense that the product of two invariants of orders $i$
and $j$ is an invariant of order $i+j$. Therefore the set of independent 
Vassiliev invariants at a given order can be divided in two subsets: those that
are products of invariants of lower orders (composite invariants), an 
those that are not (primitive invariants). We shall call the decomposition
of a Vassiliev invariant as a product of lower order Vassiliev invariants 
``factorization''. This phenomenon is most clearly exposed after 
choosing a particular kind of basis of group factors that we will call 
``canonical''. The detailed description of these bases and its significance 
to the theory of numerical knot invariants of finite type is the main goal of 
the present work.

More precisely, in this paper we shall show 
that the factorizations observed in \cite{numbers} can be resummed in a single 
exponential, which includes only the primitive Vassiliev invariants of the 
knot $C$, thus disentangling the contribution of these primitive invariants 
to all orders in perturbation theory. 
This, which is the  main result of this paper, can be regarded as an extension
of the theorems by Birman and Lin in \cite{birlin,lin,birman}
where it is proven that the coefficients of the power expansion of any
Chern-Simons or quantum group polynomial invariant is a Vassiliev 
invariant. Our main result is contained in the ``Factorization Theorem" 
presented in sect.~5 and it can be very simply stated as follows:

\vskip 0.2truein

\noindent  Let $C$ be a knot and let ${\cal H}_t^R(C,G)$ be a Chern-Simons or
quantum group polynomial invariant associated to a compact semi-simple Lie 
group $G$ and a representation $R$ of $G$ (normalized so that for the unknot 
it takes the value 1). ${\cal H}_t^R(C,G)$ is a polynomial in $t$. 
Let ${\cal W}_x^R(C,G)$ be obtained from ${\cal H}_t^R(C,G)$ by replacing the
variable $t$ by $\ex^x$ and let us consider the power series expansion of  
$\log{\cal W}_x^R(C,G)$ around $x=0$: 
\begin{equation}
\log {\cal W}_x^R(C,G) = \sum_{i=0}^\infty w_i^c x^i.
\label{maineq}
\end{equation}
Then, $w_0^c=0$ and each $w_i^c$, $i\geq 1$, is a 
{\it primitive} Vassiliev invariant relative to a canonical basis.

\vskip 0.2truein

\noindent The proof of this theorem is accomplished through the choice 
of a canonical basis for the independent group factors at each order in the 
perturbative expansion of the vacuum expectation value
of the Wilson line, $\langle\W_R(C,G)\rangle$, which is precisely, up to a
normalization factor, ${\cal H}_t^R(C,G)$, or ${\cal W}_x^R(C,G)$, above. 
The use of these  bases reveals a simplicity in the perturbative expansion that
can hardly be  grasped otherwise. 

This paper is organized as follows. Section 1 contains an elementary 
exposition of Chern-Simons quantum field theory, along with the definition of
the object of interest in this article: the Wilson loop operator. In section 3 
we introduce the general 
structure of the perturbative expansion of $\langle\W_R(C,G)\rangle$
and the definition of canonical bases. Section 4 presents some 
consequences of our having chosen a canonical basis to express the perturbative
expansion, encoded in a ``Master Equation". In Section 5 these results are 
shown to lead to the factorization of the expansion in a single exponential. 
It is noteworthy that this implies that the logarithm of the polinomial
invariant contains primitive Vassiliev Invariants only; the concept of 
primitive invariant has a simple definition in a canonical basis. Changes of
basis are analyzed in section 6. Section 7 contains our conclusions.

\section{Chern-Simons theory}
In this section we will recall a few facts on Chern-Simons gauge theory. 
Let us consider a compact semi-simple Lie group $G$ and a connection $A$ on
${\RR^3}$. The  Chern-Simons action is defined as:
\begin{equation}
S_k(A)={k\over 4\pi}\int_{\RRs^3} \tr (A\wedge dA + 
{2\over 3} A\wedge A\wedge A),
\label{action}
\end{equation}
where Tr denotes the trace in the fundamental representation of
$G$. Given a knot $C$, {\it i.e}, an embedding of $S^1$ into ${\RR^3}$, we
define the Wilson line associated to $C$ carrying a representation $R$
of $G$ as:
\begin{equation}
\W_R(C,G)=\tr_R\left[{\hbox{\rm P}} \exp \oint A \right],
\label{wilson}
\end{equation}
where ``P" stands for path ordered and the trace is to be taken 
in the representation $R$. The vacuum expectation value
is defined as the following ratio of functional integrals:
\begin{equation}
\langle \W_R(C,G) \rangle = {1\over Z_k} \int [DA] \W_R(C,G)\, {\rm e}^{iS_k(A)},
\label{vev}
\end{equation}
where $Z_k$ is the partition function:
\begin{equation}
Z_k=\int [DA]\, {\rm e}^{iS_k(A)}.
\label{parfun}
\end{equation}
The theory based on the action (\ref{action}) possesses a gauge symmetry which
has to be fixed. In addition, one has to take into account that the
theory must be regularized due to the presence of divergent integrals
when performing the perturbative expansion of (\ref{vev}). Regarding these two
problems we will follow the approach taken in \cite{numbers}. It is 
known \cite{alr} that in the Landau gauge the contribution of the one-loop 
gauge field self-energy and one-loop gauge field vertex to the perturbative 
expansion of any \vev\ can be traded by a shift in the parameter $k$ that 
multiplies the Chern-Simons action: $k\rightarrow k-C_A$, being $C_A$ the 
quadratic casimir in the adjoint representation of $G$. In so doing we do not 
need to include one loop two- or three-point gauge field subdiagrams in the 
perturbative expansion. Also, it has been shown \cite{carmelo} that 
higher-order corrections to the gauge field two- and three-point functions 
vanish.  

There is one more problem emanating from perturbative quantum field
theory which must be considered. Often, products of operators
$A_\mu(x) A_\nu(y)$ must be considered at the same point $x=y$, where
they are ambiguous. This situation can be solved \cite{witten,gmm} without 
spoiling the topological nature of the theory. In the process one needs to
introduce a framing attached to the knot which is characterized by an
integer $n$. It was shown in \cite{pert}
that working in the standard framing, $n=0$, is equivalent to ignoring 
diagrams containing collapsible propagators in the sense explained 
in \cite{pert,numbers}. 

\section{General structure of the perturbative expansion}
The facts mentioned in the previous section (exclusion of loop contributions
to gauge field two- and three-point functions, and of collapsible propagators) 
simplify considerably the perturbative analysis of the \vev\ (\ref{vev}).  As 
shown in \cite{numbers} the perturbative expansion of the vacuum
expectation value of the Wilson line (\ref{wilson}) has the form:
\begin{equation}
\langle \W_R(C,G) \rangle = d(R) \sum_{i=0}^\infty \sum_{j=1}^{d_i}
\alpha_i^j r_{ij} x^i,
\label{expansion}
\end{equation}
where $x={2\pi i \over k-C_A}$, and $d(R)$ is the dimension of the
representation $R$. The factors $\alpha_i^j$ and $r_{ij}$ in
(\ref{expansion}) incorporate all the dependence dictated from the Feynman 
rules apart from the dependence on $k$ which is contained in $x$. The power of
$x$, $i$, represents the order in perturbation theory. Of the two
factors, $r_{ij}$ and $\alpha_i^j$, the first one contains all the
group-theoretical dependence, while the second all the geometrical
dependence. The quantity $d_i$ denotes  the number of
independent group structures $r_{ij}$ which appear at order $i$. The
first values of $d_i$, $\alpha_i^j$ and $r_{ij}$ are: $\alpha_0^1=
r_{0,1}=1$, $d_0=1$ and $d_1=0$. Notice that we are shifting $k$ in the 
definition of $x$ and therefore no diagrams with
loop contributions to two and three-point functions should be
considered. In addition, there is no linear term in the expansion ($d_1=0$)
so that diagrams with collapsible propagators (isolated chords) 
should be ignored in the
sense explained in \cite{numbers}. It was proven in \cite{numbers} 
that the quantities $\alpha_i^j$ are Vassiliev invariants of order $i$.

We introduce now some vocabulary in order to classify Feynman diagrams. We
will assume that the reader is familiar with the types of trivalent Feynman
diagrams appearing in Chern-Simons perturbation theory. These diagramas are
trivalent graphs with a distinguished line called Wilson line which carry the
representation chosen. A detailed account of Chern-Simons perturbation
theory specially suited for our pourposes can be found in \cite{numbers}.  We
begin introducing the notion of ``connected" loop diagram. We will say that a
diagram is a connected loop diagram if it is  possible to go from one 
propagator (or internal line) to another without ever having to go through 
the Wilson line. If the diagram is ``disconnected'' that is not possible. In 
this second case we say that the diagram has subdiagrams, which are the 
connected components of the whole loop diagram. We say that two  subdiagrams 
are ``non-overlapping'' if we can move along the Wilson line meeting all the 
legs of one subdiagram first, and all the legs of  the other 
second. Here, ``legs'' means the propagators directly
attached to the Wilson line. If it is impossible to do that, the subdiagrams
are ``overlapping''. In fig. 1 the diagram {\it a} is connected while the
diagrams {\it b} and {\it c} are disconnected. Of these last two,
diagram {\it b} contains subdiagrams which are overlapping while 
{\it c} does not.

\begin{figure}
\centerline{\hskip.4in\epsffile{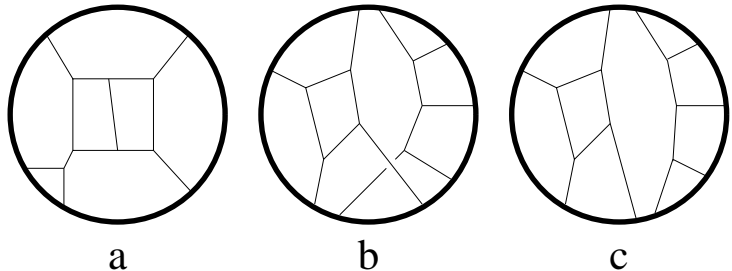}}
\caption{Examples of diagrams.}
\end{figure}

In the general expansion (\ref{expansion}) there are many possible choices of
the independent groups factors $r_{ij}$. Given all Feynman diagrams
contributing to a given order in perturbation theory some of the resulting 
group factors might be related due to the  relations among the generators
$T^a$ and the structure constants $f_{ijk}$ of semi-simple groups. From a
diagramatic point of view these relations are the so called  STU 
and IHX relations \cite{drortopo}. Since for a semi-simple group the 
structure constants can be chosen antisymmetric there is no need to 
distinguish orientation of internal three-vertices. The group factors entering 
(\ref{expansion}) are chosen to be associated to diagramas that are 
independent. Of course, many choices are possible. Each possible set of group 
factors $r_{ij}$ represents a basis. There are two simple  but far-reaching 
facts about the bases $r_{ij}$ which we summarize in two Propositions:

\vskip .5cm
{\bf Proposition 1:}
{\it It is always possible to choose a basis such that the $r_{ij}$ come from 
connected diagrams, or products of connected diagrams.} That is, if there are
subdiagrams, they can be chosen so that they do not overlap. The value of such 
an $r_{ij}$ is the product of the values of its subdiagrams. 
\vskip0.2cm
{\bf Proposition 2:}
{\it The $r_{ij}$ which are products can be chosen as products of connected 
$r_{ij}$'s of lower orders.}
\vskip .5cm

These propositions follow from two simple facts. First, using STU relations 
it is always possible to trade in a disconnected diagram overlapping 
subdiagrams by connected diagrams and disconnected diagrams containing 
non-overlapping subdiagramas. Second, if a loop diagram is not
connected, and its subdiagrams are non-overlapping,  its group factor is the
product of the group factors of its subdiagrams. This last statement follows
from the fact that if one cuts a Wilson line at a given point where no leg
is inserted the resulting matrix is a diagonal matrix.

\begin{figure}
\centerline{\hskip.4in\epsffile{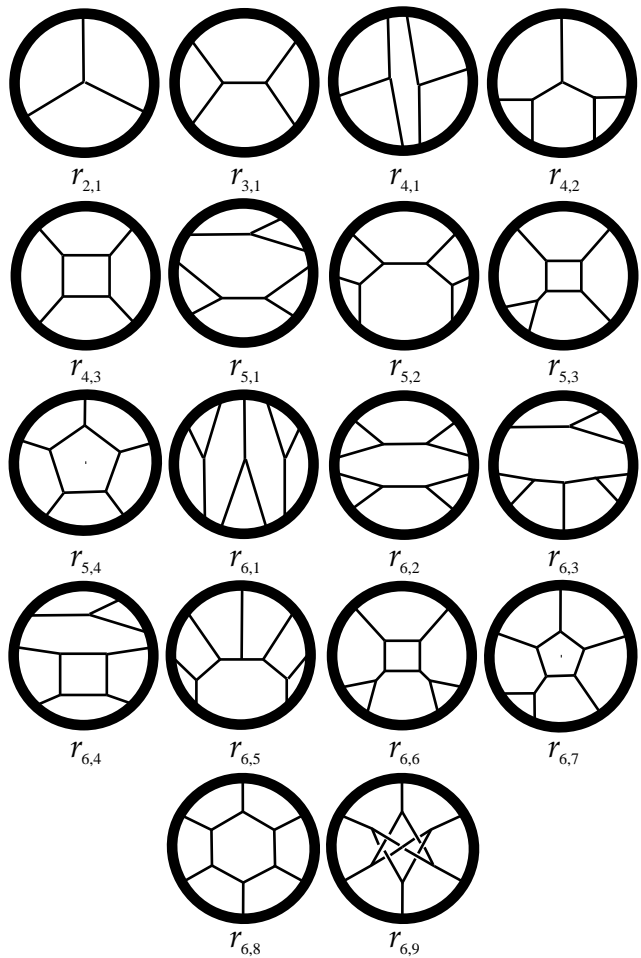}}
\caption{Example of a canonical basis up to order 6.}
\end{figure}

Propositions 1 and 2 are very important because a basis of group factors such
constructed shows the following unique feature: a given connected $r_{ij}$
begets a  whole family of other group factors at higher orders in which it
enters as  a subdiagram. A basis constructed following Propositions 1 and 2
shall be called ``a canonical basis''. The basis used in \cite{numbers}
up to order six is canonical. The diagramas chosen are reproduced in fig. 2.

\section{The Master Equation}
Let the gauge group be a product $G\otimes G'$, where $G$ and $G'$ are compact
semi-simple groups. From a path-integral  representation of the \vev s, the
following identity is obvious: \begin{equation}
\langle\W_{R\otimes R'}(C,G\otimes G')\rangle=\langle\W_R(C,G)\rangle
\langle\W_{R'}(C,G')\rangle\,\, .
\label{factorization}
\end{equation}
When combined with the choice of the same canonical basis for all the vev's, 
this equation proves to be most fruitful. In order to show this, consider an 
$r_{ij}$ composed of several connected subdiagrams which we denote by 
$r_{ij}^{(p)}$ with $p=1,\ldots,\#(ij)$. Some of them may be identical. In a 
canonical basis these 
$r_{ij}^{(p)}$ are elements of the basis at lower orders, and therefore are 
associated with geometrical factors which we denote by $(\alpha_i^j)^{(p)}$. 
If the Lie group is simple, we know that
\begin{equation}
r_{ij}(G)=\prod_{p=1}^{\#(ij)} r_{ij}^{(p)}(G) ,
\label{erres1}
\end{equation}
but if the Lie group is a product, we find that
\begin{equation}
r_{ij}(G\otimes G')=\prod_{p=1}^{\#(ij)}\left( r_{ij}^{(p)}(G)+ r_{ij}^{(p)}(G')
\right).
\label{erres2}
\end{equation}
Inserting eq.~(\ref{expansion}) in eq.~(\ref{factorization}) and putting 
things together, we arrive at the ``Master Equation'':
\begin{eqnarray}
& &\sum_{i=0}^{\infty}\sum_{j=1}^{d_i}\alpha_i^j(C)
\prod_{p=1}^{\#(ij)}\left( r_{ij}^{(p)}(G)x^i + 
r_{ij}^{(p)}(G'){x'}^i \right)=\nonumber \\
& &\left( \sum_{k=0}^{\infty}\sum_{l=1}^{d_k}
\alpha_k^l(C) r_{kl}(G)x^k\right)\left(\sum_{m=0}^{\infty}\sum_{n=1}^{d_m}
\alpha_m^n(C) r_{mn}(G'){x'}^m\right) \,\, ,
\label{master}
\end{eqnarray}
where $x=2\pi i/(k-C_A)$ and $x'=2\pi i/(k-C_A')$.
In eq. (\ref{erres1}), (\ref{erres2}) and (\ref{master}) we have shown
 explicitly the fact
that the group factors $r_{ij}$ depend only on the group-theoretical data while
the geometrical factors $\alpha_i^j$ depend only on the knot $C$. Actually, we
would ought to have indicated the representations $R$ and $R'$ in the group
factors. We did not do it to avoid  a cumbersome notation but it
certainly should be understood.
The matching of the two 
polynomials in $x$ and $x'$ in (\ref{master}) produces an infinite 
string of identities relating $\alpha_i^j$'s at a given order with products
of the $(\alpha_k^l)^{(p)}$ of its components. The general result is as 
follows. Let us consider a composite $r_{ij}$ that consists of $p_{ij}^{(k)}$ 
connected non-overlapping subdiagrams of some type $k$, with $k=1, \ldots,N$. 
This means in particular that
\begin{equation}
\sum_{k=1}^N p_{ij}^{(k)} = \#(ij)\,\, .
\label{check}
\end{equation}
(The only purpose of this formula is to clarify the notation). We call 
$r_{ij;k}$ the connected subdiagram of type $k$, which
in a canonical basis is also an element of the basis at a lower order, and 
therefore is associated to a geometrical factor denoted by
$\alpha_{i\,\,;k}^j$. In other words, the element $r_{ij}$ contains the
following connected subdiagrams:
\begin{eqnarray}
r_{ij}^{(1)}&=&r_{ij}^{(2)}=\ldots=r_{ij}^{(p_{ij}^{(1)})}\equiv r_{ij;1}\, ,
\nonumber \\
r_{ij}^{(p_{ij}^{(1)}+1)}&=&r_{ij}^{(p_{ij}^{(1)}+2)}=\ldots=
r_{ij}^{(p_{ij}^{(1)}
+p_{ij}^{(2)})}\equiv r_{ij;2}\, , \\ &{\rm etc.}& \nonumber
\label{explain}
\end{eqnarray}
The Master Equation (\ref{master}) 
leads to a formula for the $\alpha_i^j$ associated with 
our composite $r_{ij}$:
\vskip .5cm
{\bf THEOREM 1}
\vskip0.2cm
\begin{equation}
\alpha_i^j= \prod_{k=1}^N {1\over p_{ij}^{(k)}!}\,\left(\alpha_{i\,\,;k}^j
\right)^{p_{ij}^{(k)}}\,\, .
\label{theor1}
\end{equation}
Note that this is true only for a canonical choice of basis. This general 
result is the key to the next sections.  

\section{Factorizations}
When the approach described in the preceding sections was first proposed in
\cite{numbers}, the canonical basis did not include any elements at order $1$. 
The reason is that the only contribution to the \vev\ at that order would be 
the framing factor, which is not intrinsic to the knot. In that paper the 
interest was centered in extracting numerical knot invariants from the 
perturbative expansion, so it was natural to exclude the framing factor 
and its corresponding element of the basis. 

We can include the framing factor if we enlarge our basis with $r_{11}=C_2$, 
the quadratic casimir. The geometrical factor associated with this new group 
structure will be denoted by $n$ for it is the framing. We are including these 
new elements here in order to show the simplest application of 
eq.~(\ref{theor1}).

According to Propositions 1 and 2 there is a new basis in which each $r_{kl}$ 
originates a family of elements of the form $r_{kl}C_2$, $r_{kl}C_2^2$ and 
so on. This set of elements of the new basis originated by $r_{kl}$ can be 
called ``the $kl$-$C_2$-family''. Let us focus on a given $kl$-$C_2$-family 
and consider its contribution to $\langle\W_R(C,G)\rangle$. Theorem 1 shows that 
this contribution is
\begin{equation}
d(R)\left( \alpha_k^l\,r_{kl}\,x^k +  \alpha_k^l\,r_{kl}\,C_2\,n\,x^{k+1} +
 \alpha_k^l\,r_{kl}\,C_2^2\,{1\over 2}n^2\,x^{k+2} + \ldots \right),
\label{contri}
\end{equation}
and the following terms are the expansion of an exponential. It follows that 
eq.~(\ref{contri}) is equal to:
\begin{equation}
d(R) \alpha_k^l\,r_{kl}\,x^k\, \ex^{C_2 n x}.
\label{exp1}
\end{equation}
Therefore, repeating the same argument for each $kl$-$C_2$-family (which are 
all independent because the $r_{kl}$ are) we get:
\vskip .5cm
{\bf THEOREM 2}
\vskip0.2cm
\begin{equation}
\langle\W_R(C,G)\rangle=\langle\W_R(C,G)\rangle\big|_{n=0}\, \ex^{C_2 n x}.
\label{exp2}
\end{equation}
This agrees exactly with the non-perturbative result \cite{witten} and can be
regarded as a simplified proof of the factorization of the framing factor shown
in \cite{pert}. Now, we can use a similar approach to factorize more structures.
The idea is the same: the selection of a  given connected element of our
canonical basis, the  ``dressing'' of another  element of the basis with copies
of our selected  subdiagram (thus creating a family of diagrams in the sense
explained above)  and the repeated use of Theorem 1. The rest of this section
dwells on this  subject. 

The reasoning that led to the factorization of the $C_2$ structure can be 
applied without changes to any other connected $r_{ij}$. The only peculiarity
of $C_2$ is that its addition to an element of the basis at order $k$, say
$r_{kl}$, leads to another element at order $k+1$, its ``first descendant'' in
the $kl$-$C_2$-family. If we want to factorize $r_{21}$, the members of the 
$kl$-$r_{21}$-family would have orders $k+2$, $k+4$ and so on in double 
steps. Proposition 2 says that our basis can be constructed so that it contains 
this new family. 

Let $r_{k+2q, l}$ be a member of the $kl$-$r_{21}$-family generated by 
$r_{kl}$ and a $q$-fold insertion of $r_{21}$, {\it i.e.}
\begin{equation}
r_{k+2q, l} = r_{kl}\,\left(r_{21}\right)^q\, .
\label{erre3}
\end{equation}

Theorem 1, and the fact that we have chosen a canonical basis, enables us to 
prove that the contribution of this family to $\langle\W_R(C,G)\rangle$ is   
\begin{eqnarray}
& &d(R)\,\left(\alpha_k^l\,r_{kl}\,x^k+\alpha_k^l\,\alpha_2^1\,r_{kl}\,
r_{21}\,x^{k+2}+{1\over 2}\alpha_k^l\,(\alpha_2^1)^2\,r_{kl}\,\,
(r_{21})^2\,x^{k+4}+\ldots \right)\nonumber\\
& &=d(R)\,\alpha_{kl}\,r_{kl}\,x^k\,\ex^{\alpha_2^1\,r_{21}\,x^2}.
\label{facalpha21}
\end{eqnarray}

We can repeat the same steps for all $kl$-$r_{21}$-families because they are 
all independent; in so doing we arrive at:
\vskip .5cm
{\bf THEOREM 3}
\vskip0.2cm
\begin{equation}
\langle\W_R(C,G)\rangle=\langle\W_R(C,G)\rangle\big|_{\alpha_2^1=0}\, 
\ex^{\alpha_2^1\,r_{21}\, x^2},
\label{theo3}
\end{equation}
and nothing hinders the generalization of this theorem to what can be seen as
the full expression of the concept of perturbative factorization of the \vev :
\vskip .5cm
{\bf FACTORIZATION THEOREM}
\vskip0.2cm
\begin{equation}
\langle\W_R(C,G)\rangle=d(R)\,\exp\left\{\sum_{i=1}^{\infty}\sum_{j=1}^
{\hat{d}_i}\alpha_i^{j\,\,c}(C)\,r_{ij}^c(G)\,x^i\right\},
\label{theofacto}
\end{equation} 
where $r_{ij}^c$ denotes the connected elements of the basis, and 
$\alpha_i^{j\,\,c}$ their corresponding geometrical factors. These 
$\alpha_i^{j\,\,c}$ do {\it not} correspond uniquely to connected diagrams 
because the $r_{ij}^c$ include geometric factors from both connected and 
disconnected Feynman diagrams. The symbol $\hat{d}_i$ stands for the number of 
connected elements in the canonical basis at order $i$. 

Equation~(\ref{theofacto}) can be written in the form
\begin{equation}
\log\left( {1\over d(R)}\langle\W_R(C,G)\rangle\right)=\sum_{i=1}^{\infty}
\sum_{j=1}^{\hat{d}_i}\alpha_i^{j\,\,c}(C)\,r_{ij}^c(G)\,x^i,
\label{logarithm}
\end{equation}
which is the result announced in (\ref{maineq}). Equation (\ref{logarithm})
is reminiscent  of the 
well known fact in quantum field theory that the logarithm of the generating 
functional can be expanded in terms of connected diagrams only. The relevance 
of this formula to the theory of knot invariants comes from the identification
of the \vev\ of the Wilson line as a polynomial invariant. Actually we are in 
a much more general situation because we are considering an arbitrary 
semi-simple gauge group, so our vev is in some sense the most general 
polynomial invariant possible. Therefore, eqs.~(\ref{theofacto}) or 
(\ref{logarithm}) prove that if a canonical basis is chosen, the logarithm of 
a polynomial knot invariant can be expanded in terms of the primitive 
Vassiliev invariants of the knot only. 

The primitive Vassiliev invariants $\alpha_i^{j\,\,c}(C)$ have been computed
up to order six for all prime knots up to six crossings \cite{numbers}
and for arbitrary torus knots \cite{torus}. It was conjectured in 
\cite{numbers} that there exist a normalization for the $\alpha_i^{j\,\,c}(C)$ 
in which they are integer-valued. The integral expressions for
$\alpha_2^1(C)$ and $\alpha_3^1(C)$ were first presented in \cite{gmm}
and in \cite{numbers} respectively. Properties of these two primitive Vassiliev
invariants have been studied in \cite{drortesis} and in  \cite{alemanes}.

\section{Change of basis}
In this section we want to investigate to what extent the previous results are
independent of the basis chosen. We are thus led to consider changes of 
basis. First we treat changes of canonical basis. Let $B$ and $B'$ be two 
different canonical bases, being $r_{ij}$ and $r_{ij}'$ its elements. The most 
general change of canonical basis is
\begin{equation}
r_{ij}'= N_j^k\,r_{ik},
\label{changer}
\end{equation}
where $N$ is a $d_i \times d_i$ matrix \footnote {There is an $N$ at each order
$i$, but we are not indicating this fact.} yet to be determined. We are 
assuming that the vectors in this space are written as columns. To see that 
(\ref{changer}) is the most general change between $B$ and $B'$, consider a 
possible extra term in the right hand side. It has to be a product of several 
$r_{kl}\in B$ such that the sum of the orders of its factors is $i$. But by 
definition of canonical basis such a product is also an element of $B$ at 
order $i$ and therefore the extra term can be absorbed into the first 
term. This shows that eq.~(\ref{changer}) is indeed the most general change of 
canonical basis.

It is elementary to prove that, in order to preserve the independence of the
elements of $B$ at each order, the matrix $N$ has to be non-singular:
\begin{equation}
\det N\neq 0\,\,.
\label{nonsing}
\end{equation}
We now analyze the effect
of the change of basis (\ref{changer}) on the geometric factors $\alpha_i^j$. 
The vectors in this space of numerical factors are written as rows. The 
starting point is the invariance of the \vev\ under changes of 
basis. Therefore, at each order $i$ it is true that,
\begin{equation}
\sum_{j=1}^{d_i}\,\alpha_i^j\,r_{ij}=
\sum_{j=1}^{d_i}\,{\alpha'}_i^j\,{r'}_{ij}\, ,
\label{sea}
\end{equation}
and it follows that the $\alpha_i^j$ transform ``contravariantly'':
\begin{equation}
{\alpha'}_i^j=\,\alpha_i^k\,\left(N^{-1}\right)_k^j\,\,.
\label{changealpha}
\end{equation}
In general, at order $i$ the change of canonical basis involves a matrix 
$N\in GL(d_i,\QQ)$. We do not know how to characterize these matrices
completely.  Only a small subset of the whole linear group is relevant. For
example we can  discard $N$'s which are mere permutations, or diagonal, since
they do not  lead to essentially new bases. More important, the elements of the
basis $B'$ will be linear combinations of the elements of $B$, but these linear 
combinations must be interpretable as {\it new diagrams\/} because we want
$B'$ to be canonical. In other words, 
we would have to investigate which linear combinations of independent diagrams
can be written as a single diagram. The resulting subgroup of $GL(d_i,\QQ)$
would be equivalent to the space of independent canonical bases at order $i$.

The $N$'s have some properties that do not depend on these details. Let us 
order the elements of a canonical basis $B$ at order $i$ as follows:
\begin{equation}
B_i=\left\{r^c_{i1}, \ldots ,r^c_{i, \hat{d}_i}, r_{i, \hat{d}_i+1},\ldots, 
r_{i, d_i}\right\}\, , 
\label{order}
\end{equation}
{\it i.e.\/} the connected elements before, and the disconnected elements 
after. A given disconnected element of $B$ can be written as a product of 
elements of $B$ of lower orders, 
\begin{equation}
r_{ij}=r_{kl}\,r_{i-k, s}\, ,
\label{decomp}
\end{equation}
where $l$ and $s$ depend implicitly on $j$, but this will not be relevant in 
what follows. If we write the identity (\ref{decomp}) in a new canonical basis 
$B'$, it reads
\begin{equation}
N_j^p\,{r'}_{ip}=N_l^q\,{r'}_{kq}\,N_s^t\,{r'}_{i-k, t}\, .
\label{decompprima}
\end{equation}
Note that the $N$'s operate on different spaces. What we
have on the right hand side is a linear combination of elements of $B'$
at order $i$, because all these products of $r_{kq}$'s times $r_{i-k, t}$'s
must be elements of $B'$ (it is canonical). On the left hand side we have other
elements of the same basis at the same order. Therefore eq.~(\ref{decompprima})
is a contradiction unless the $r_{ip}'$ are themselves products and, therefore,
correspond to disconnected diagrams. The conclusion is that in a change of 
canonical basis, the disconnected $r$'s in $B$ come from the disconnected 
$r$'s in $B'$.

This result can be summarized in the following symbolic representation of a 
matrix $N$ valid for a change of canonical basis:
\begin{equation}
N=\left( \begin{tabular}{c|c}
		C $\rightarrow$ C & NC $\rightarrow$ C \\ \hline
		0 & NC $\rightarrow$ NC
	 \end{tabular} \right)
\label{matriz}
\end{equation} 
where $C$ and $NC$ mean ``connected'' and ``non-connected'' respectively. In 
this notation a change of basis would be written as
\begin{equation}
\left(\begin{array}{c} r^c \\ r^{nc} \end{array}\right) =\,
\left(\begin{tabular}{c|c}
		 A & B \\ \hline
		C & D
	 \end{tabular} \right)\,
\left(\begin{array}{c} {r'}^c \\ {r'}^{nc} \end{array} \right)
\label{nota}
\end{equation}
where $C=0$. The blocks in the diagonal are square matrices because all 
canonical bases must have the same number $\hat{d}_i$ of connected elements 
at a given order $i$. These matrices have an interesting property: they form a 
subgroup of $GL(d_i,\QQ)$. The inverse of a given element is 
\begin{equation}
\left( \begin{tabular}{c|c}
		 A & B \\ \hline
		0 & D
	 \end{tabular} \right)^{-1}=
\left( \begin{tabular}{c|c}
		${\rm A}^{-1}$ & $-{\rm A}^{-1}{\rm B}{\rm D}^{-1}$ \\ \hline
		$0$ & ${\rm D}^{-1}$
	 \end{tabular} \right)
\label{inverse}
\end{equation} 
and the determinant is
\begin{equation}
\det \left( \begin{tabular}{c|c}
		 A & B \\ \hline
		 0 & D
	 \end{tabular} \right)= \det{\rm A}\,\det{\rm D}\,\,.
\label{det}
\end{equation}
Therefore, the matrix is non-singular if and only if its diagonal blocks are
non-singular. We are assuming that the matrix represents a valid change of 
basis, thus nothing is singular and the inverses in eq.~(\ref{inverse}) do
exist.

\subsection{Diagrammatic interpretation}
We can sharpen the previous result by analyzing the constraints that the 
diagrammatic origin of the $r$'s imposes on their algebra. For example, the 
product of two $r$'s is interpretable as a diagram in an obvious way 
(actually, as an equivalence class of diagrams). The big constraint is on the 
sum of two $r$'s. We are now investigating the algebra of diagrams which are 
either connected or disconnected and non-overlapping, independently of their 
being independent or not, {\it i.e.\/} of being elements of a canonical 
basis or not. We shall call these diagrams generally $r$, and will not consider
diagrams with overlapping subdiagrams at all. 

The question is: when is the sum of two $r$'s interpretable as a diagram? The 
answer comes in two parts:
\vskip .5cm
{\bf Lemma}
\begin{enumerate}
\item Let $r_1$ and $r_2$ be two different connected diagrams of the same 
order $i$. Their sum or difference $r_1 \pm r_2$ exists as a diagram of order 
$i$ if and only if $r_1$ and $r_2$ are two of the terms in an STU or IHX 
relation. The sign in $\pm$ depends on which two terms of those relations
are $r_1$ and $r_2$ related to. 
\item If $r_1$ and $r_2$ are not connected, they can differ only in a single
subdiagram. Again, the subdiagrams that are different must be two of the terms
in an STU or IHX relation.
\end{enumerate}

More complicated linear combinations of $r$'s that are interpretable as 
diagrams can always be decomposed in elementary steps, each of them satisfying
the Lemma. We introduce now a new concept. A linear combination of $r$'s that 
can be interpreted as a diagram shall be called a ``valid'' linear 
combination. To be completely rigurous, a valid linear combination should be 
written in an unambiguous way by using parentheses to indicate which $r$'s 
are added to which other $r$'s and in what order. For example, 
\begin{equation}
r_1+r_2-r_3 = (r_1+r_2)-r_3\quad {\rm or}\quad r_1+(r_2-r_3)\quad {\rm or} 
\quad(r_1-r_3)+r_2 \,\,.
\label{minor}
\end{equation}
Which one of the three possibilities is the good one depends on the diagrams. 
In this sense, the addition of $r$'s is commutative but not associative. This 
may be a minor point, because the {\it numerical values} of the $r$'s can be
added freely, but we are now focusing on the formal properties of the $r$'s
as elements of an ``algebra'' \footnote{Strictly speaking they do not form an 
algebra, whence the quotation marks.}. In general, when we write 
a valid linear combination of $r$'s 
we shall assume that there is an ordering of the additions such that each 
step complies with the Lemma. This ordering depends on the particular diagrams 
in the linear combination and thus will not be indicated in general. 

\subsection{Arithmetic of diagrams}
The previous subsection establishes the rules for an arithmetic of diagrams. 
The $r$'s are elements of an algebraic structure in which we
can always multiply, but not always add or substract (see the Lemma above). As 
for the division, an expression like $r_1/r_2$ is an $r$ only if $r_2$ is a 
subdiagram of $r_1$, in which case we say that $r_2$ divides $r_1$. There is a 
neutral element for the product: the empty diagram. Therefore no diagrammatic 
interpretation exists for $1/r_2$.

The set of $r$'s with the addition is not even a group, for not all diagrams 
can be added. This lack of structure precludes an abstract formulation
of the algebra of group factors $r$. We need to know which is the
diagram that a given $r$ represents in order to ascertain if it can be added
to another $r$ or not. Nevertheless we can prove some theorems which have a
bearing in our investigation. 
\vskip .5cm
{\bf THEOREM 4}
\vskip0.2cm
Let $r_1$ and $r_2$ be two diagrams of the same order $i$ that can be 
added or substracted, and $r_3$, $r_4$ be two diagrams whose orders add to $i$.
Then,  
\begin{equation}
r_1 \pm r_2 = r_3\,r_4 \,\,\Longleftrightarrow\,\, \left(r_3 \big| r_1 \quad 
{\rm and}\quad r_3 \big| r_2 \right)
\quad {\rm or}\quad \left(r_4 \big| r_1 \quad{\rm and}\quad r_4 \big| r_2
\right)
\label{chorras}
\end{equation}
where $a|b$ means that $a$ divides $b$. The proof follows from the 
Lemma. An immediate generalization of this theorem is that a valid linear
combination of $r$'s is disconnected non-overlapping if and only if so is
each $r$ in the linear combination. A similar conclusion holds for valid
linear combinations of connected diagrams. We want to gather these important
conclusions in three remarks:
\begin{enumerate}
\item No valid linear combination of $r$'s includes connected
and disconnected non-overlapping diagrams at the same time.
\item A valid linear combination of connected diagrams is a 
connected diagram. 
\item A valid linear combination of disconnected 
non-overlapping diagrams is a disconnected non-overlapping diagram. All of them
have the same number of subdiagrams.
\end{enumerate}

We can say that ``a valid linear combination of diagrams conserves the number 
of components''. A closely related result is that a valid change of canonical 
basis must be represented by a block-diagonal matrix:
\vskip .5cm
{\bf THEOREM 5}
\vskip0.2cm
Let $N$ be a matrix corresponding to a change of canonical basis at order $i$, 
written in the form
\begin{equation}
N=\,\left( \begin{tabular}{c|c}
		 A & B \\ \hline
		C & D
	 \end{tabular} \right)
\label{otraN}
\end{equation}
where A is a non-singular $\hat{d}_i\times \hat{d}_i$ matrix, and D a 
non-singular $(d_i-\hat{d}_i)\times (d_i-\hat{d}_i)$ matrix. Then, 
\begin{equation}
{\rm B}={\rm C}=0 \,\,.
\label{ceros}
\end{equation}
To prove this theorem notice that under a change of basis the connected $r$'s
transform as (see  eq.~(\ref{nota}) for notation)
\begin{equation}
r^c={\rm A}{r'}^c+{\rm B}{r'}^{nc}\,\,.
\label{fedup}
\end{equation}
Given that the l.h.s. is a diagram we observe that the r.h.s. is a valid linear
combination of diagrams. But the Lemma implies that all diagrams in the r.h.s.
must be connected, or all of them disconnected non-overlapping, {\it i.e.\/}
either A=0 or B=0. It is clear from the Lemma that the only option is
B=0. A similar argument establishes that C=0, which we already knew from
previous considerations, see eq.~(\ref{matriz}).

The final picture for a valid $N$ is
\begin{equation}
N=\, \left( \begin{tabular}{c|c}
		A & 0 \\ \hline
		0 & D
	 \end{tabular} \right)
\label{finalN}
\end{equation}
with A and D non-singular square matrices. In words, the connected $r$'s
transform independently from the disconnected $r$'s; they never mix in a
change of canonical basis. 

In particular this shows that eq.~(\ref{theofacto}) is consistent under a 
change of canonical basis. No matter what canonical basis we choose, the only 
relevant diagrams are those that contribute to the $r^c_{ij}$. As for the 
change from a canonical $B$ to a non-canonical $B'$, we have little to 
say. The concept of factorization as described here only makes sense for 
canonical bases, and only in this case we have eq.~(\ref{theofacto}). 

We regard canonical bases as privileged systems of reference in which the
perturbative expansion is at its simplest.

\section{Conclusions}
We have shown that the perturbative expansion of the vev of a Wilson line in 
Chern-Simons quantum field theory shows a striking simplicity if presented in
terms of a canonical basis for the group factors. These bases provide a simple
characterization  of primitive Vassiliev invariants: they are the geometric 
factors associated to the connected elements of a canonical basis. 

Within this framework it is possible to factorize the whole perturbative 
expansion in separate contributions from each primitive invariant. Each of
these factors can be resummed in an exponential. The relevance of this result 
for the theory of numerical knot invariants is that the logarithm of a 
polynomial invariant contains only the primitive invariants. 

This work opens a variety of investigations. One should study turther the
algebraic structure of the set of canonical bases. As follows from the
discussion in sect. 6, at each order $i$ this set is characterized by a
subgroup of $GL(\hat d_i,\QQ)$. Methods to find the groups corresponding to
each order should be investigated.

Another important extension of our work  consists of the study of
factorization in the context of $n$-component links. Vassiliev invariants for
$n$-component links have not been much studied and it is very likely that some
of the ideas behind factorization can be used  in the
organization of their structure. The extension of the concept of canonical basis
and primitiveness should be explored in that case following a similar analysis to
the one presented in this paper. Work on this and other aspects of Vassiliev
invariants associated to $n$-component links will be reported elsewhere.

\vskip 0.5truein
\begin{center}
{\bf Acknowledgements}
\end{center}
M.A. is indebted to Prof. John Negele and the CTP for hospitality, and to the
Spanish CICYT for financial support. J.M.F.L. thanks E. P\'erez for many
intersting discussions on Vassiliev invariants.

\end{document}